\documentclass[12pt,reqno]{amsart}
\usepackage{mathrsfs}
\usepackage{amsfonts}
\usepackage{latexsym,amsmath,amsfonts,amssymb,amsthm}
 \theoremstyle{plain}
\newtheorem{Thm}{Theorem}
\newtheorem{Lem}{Lemma}

\theoremstyle{definition}
\newtheorem*{Ack}{Acknowledgment}
\theoremstyle{remark}
\newtheorem*{Rem}{Remark}

\def\Z{\mathbb Z}
\def\N{\mathbb N}

\def\S{\mathcal{S}}
\def\T{\mathcal{T}}
\def\cN{\mathcal{N}}

\def\1{{\bf 1}}
\def\pmod #1{\ ({\rm mod}\ #1)}

\begin{document}

\title{On the integers not of the form $p+2^a+2^b$}
\author{Hao Pan}
\email{haopan79@yahoo.com.cn}
\address{Department of Mathematics, Nanjing University,
Nanjing 210093, People's Republic of China}  \maketitle
\begin{abstract}
We prove that
$$
|\{1\leqslant n\leqslant x:\,n\text{ is odd and not of the form
}p+2^a+2^b\}|\gg x\cdot\exp\bigg(-C\log x\cdot\frac{\log\log\log\log
x}{\log\log\log x}\bigg),
$$
where $C>0$ is an absolute constant.
\end{abstract}

\section{Introduction}
\setcounter{equation}{0} \setcounter{Thm}{0} \setcounter{Lem}{0}
\setcounter{Cor}{0}

As early as 1849, Polignac conjectured that every odd integer
greater than 3 is the sum of a prime and a power of 2. Of course,
Polignac's conjecture is not true, since 127 is an evident
counterexample. In 1934, Romanoff \cite{Romanoff34} proved that the
sumset
$$
\{p+2^b:\, p\text{ is prime},\ b\in\N\}
$$
has a positive lower density. And in the other direction, van der
Corput \cite{vanderCorput50} proved that the set
$$
\{n\geqslant 1:\,n\text{ is odd and not of the form }p+2^b\}
$$
also has a positive lower density. In fact, with help of a covering
system, Erd\H os \cite{Erdos50} found that every positive integer
$n$ with $n\equiv 7629217\pmod{11184810}$ is not representable as
the sum of a prime and a power of 2.

In \cite{Crocker71}, Crocker proved that there exist infinitely many
odd positive integers $x$ not of the form $p+2^a+a^b$. One key of
Crocker's proof is the following observation: If $b-a=2^st$ with
$s\geqslant 0$ and $2\nmid t$, then
$2^a+2^b\equiv0\pmod{2^{2^s}+1}$. And Crocker also constructed a
suitable covering system to deal with the case $a=b$. In
\cite{SunLe01}, Sun and Le discussed the integers not of the form
$p^\alpha+c(2^a+2^b)$. And subsequently, Yuan \cite{Yuan04} proved
the there exist infinitely many positive odd integers $x$ not of the
form $p^\alpha+c(2^a+2^b)$.

Let
$$
\cN=\{n\geqslant 1:\,n\text{ is odd and not of the form }p+2^a+2^b\}
$$
and
$$
\cN_*=\{n\geqslant 1:\,n\text{ is odd and not of the form
}p^\alpha+2^a+2^b\}.
$$
Erd\H os asked whether $|\cN\cap[1,x]|\gg x^\epsilon$ for some
$\epsilon>0$. And as Granville and Soundararajan
\cite{GranvilleSoundararajan98} mentioned, this is true under the
assumption that there exist infinitely many $m_1<m_2<m_3<\ldots$
satisfying all $2^{2^{m_i}}+1$ are composite and $\{m_{i+1}-m_i\}$
is bounded. Erd\H os even suggested \cite[A19]{Guy09} that
$|\cN\cap[1,x]|\geqslant Cx$ for a constant $C>0$, though it seems
that the covering congruences could not help here. In
\cite{ChenFengTemplier08}, Chen, Feng and Templier proved that
$$
\limsup_{x\to\infty}\frac{|\cN_*\cap[1,x]|}{x^{1/4}}=+\infty
$$
if there exist infinitely many $m$ satisfying $2^{2^m}+1$ is
composite, and
$$
\limsup_{x\to\infty}\frac{|\cN_*\cap[1,x]|}{\sqrt{x}}>0
$$
if there are only finite many $m$ satisfying $2^{2^m}+1$ is prime.
Recently, in his answer to a conjecture of Sun, Poonen
\cite{Poonen09} gave a heuristic argument which suggests that for
each odd $k>0$,
$$
|\{1\leqslant n\leqslant x:\,n\text{ is odd and not of the form
}p+2^a+k\cdot2^b\}|\gg_\epsilon x^{1-\epsilon}
$$
for any $\epsilon>0$, where $\gg_\epsilon$ means the implied
constant only depends on $\epsilon$.

On the other hand, using Selberg's sieve method, Tao \cite{Tao}
proved that for any $\mathcal {K}\geqslant 1$ and sufficiently large
$x$, the number of primes $p\leqslant x$ such that $|kp\pm ja^i|$ is
composite for all $1\leqslant a,j,k\leqslant \mathcal {K}$ and
$1\leqslant i\leqslant \mathcal {K}\log x$, is at least $C_{\mathcal
{K}} x/{\log x}$, where $C_{\mathcal {K}}$ is a constant only
depending on $\mathcal {K}$. Motivated by Tao's idea, in this short
note, we shall unconditionally prove that
\begin{Thm}
\label{t1} $$ |\cN_*\cap[1,x]|\gg x\cdot\exp\bigg(-C\log
x\cdot\frac{\log\log\log\log x}{\log\log\log x}\bigg),
$$
where $C>0$ is an absolute constant.
\end{Thm}

Clearly Theorem \ref{t1} implies $|\cN_*\cap[1,x]|\gg_\epsilon
x^{1-\epsilon}$ for any $\epsilon>0$. The proof of Theorem \ref{t1}
will be given in the next section. And unless indicated otherwise,
the constants implied by $\ll$, $\gg$ and $O(\cdot)$ are always
absolute.

\section{Proof of Theorem \ref{t1}}
\setcounter{equation}{0} \setcounter{Thm}{0} \setcounter{Lem}{0}
\setcounter{Cor}{0}

Since
$$
|\{1\leqslant n\leqslant x:\, n\text{ is of the form
}p^\alpha+2^a+2^b\text{ with }\alpha\geqslant 2\}|=O(\sqrt{x}(\log
x)^3),
$$
we only need to show that
\begin{align*}
|\cN\cap[1,x]|\gg x\cdot\exp\bigg(-C\log
x\cdot\frac{\log\log\log\log x}{\log\log\log x}\bigg).
\end{align*}

The following two lemmas are easy applications of the Selberg sieve
method (cf. \cite[Theorems 3.2 and 4.1]{HalberstamRichert74},
\cite[Theorem 7.1]{Nathanson96}).
\begin{Lem}\label{l1} Suppose that $W\geqslant 1$ and $\beta$ are integers with $(\beta,W)=1$. Then
$$
|\{1\leqslant n\leqslant x:\,Wn+\beta\text{ is
prime}\}|\leqslant\frac{C_1x}{\log x}\prod_{p\mid
W}\bigg(1-\frac{1}{p}\bigg)^{-1},
$$
where $C_1$ is an absolute constant.
\end{Lem}
\begin{Lem}\label{l2} Suppose that $x$ is a sufficiently large
integer. Suppose that $p_1,p_2,\ldots,p_h$ are distinct primes less
than $x^{\frac18}$. Then
$$
|\{1\leqslant n\leqslant x:\,n\not\equiv0\pmod{p_j}\text{ for every
}1\leqslant j\leqslant h\}|\leqslant
C_2x\prod_{j=1}^h\bigg(1-\frac{1}{p_j}\bigg),
$$
where $C_2$ is an absolute constant.
\end{Lem}
The following lemma is due to Ford, Luca and Shparlinski
{\cite[Theorem 1]{FordLucaShparlinski09}}.
\begin{Lem}\label{l3}
The series
$$
\sum_{n=1}^\infty\frac{(\log n)^\gamma}{P(2^n-1)}
$$
converges for any $\gamma<1/2$, where $P(n)$ denotes the largest
prime factor of $n$.
\end{Lem}
Let
$$
C_3=\sum_{p\text{ prime}}\frac{1}{P(2^p-1)}.
$$
Suppose that $x$ is sufficiently large. Let
$$
K=\left\lfloor\frac{\log\log\log x}{100\log\log\log\log
x}\right\rfloor
$$
and $L=\log(2^9C_1C_2K)+2C_3$, where $\lfloor
\theta\rfloor=\max\{z\in\Z:\, z\leqslant \theta\}$.

Let $u=e^{e^{K(L+1)}}$. By the Mertens theorem (cf. \cite[Theorem
6.7]{Nathanson96}), we know that
$$
\sum_{\substack{p\leqslant u\\ p\text{ prime}}}\frac{1}{p}=\log\log
u+B+O\bigg(\frac{1}{\log u}\bigg)=K(L+1)+O(1).
$$
where $B=0.2614972\ldots$ is a constant. So we may choose some
distinct odd primes less than $u$
$$p_{1,1},\ldots,p_{1,h_1}; p_{2,1},\ldots,p_{2,h_2};\ldots;
p_{K,1},\ldots,p_{K,h_K}$$ such that
$$
\sum_{j=1}^{h_i}\frac{1}{p_{i,j}}\geqslant L
$$
for $1\leqslant i\leqslant K$. Let $q_{i,j}=P(2^{p_{i,j}}-1)$ for
$1\leqslant i\leqslant K$ and $1\leqslant j\leqslant h_i$. Clearly
these $q_{i,j}$ are all distinct. Now,
$$
\sum_{j=1}^{h_i}\log\bigg(1-\frac{1}{p_{i,j}}\bigg)\leqslant-\sum_{j=1}^{h_i}\frac{1}{p_{i,j}},
$$
whence
$$
\prod_{j=1}^{h_i}\bigg(1-\frac{1}{p_{i,j}}\bigg)\leqslant e^{-L}.
$$
And
$$
\prod_{i=1}^K\prod_{j=1}^{h_i}\bigg(1-\frac{1}{q_{i,j}}\bigg)^{-1}\leqslant\prod_{i=1}^K\prod_{j=1}^{h_i}\bigg(1+\frac{2}{q_{i,j}}\bigg)\leqslant\bigg(\frac{\sum_{i=1}^K\sum_{j=1}^{h_i}(1+2/q_{i,j})}{\sum_{i=1}^Kh_i}\bigg)^{\sum_{i=1}^Kh_i}\leqslant
e^{2C_3}.
$$

Let
$$
W_{1,i}=\prod_{j=1}^{h_i}q_{i,j}
$$
for $1\leqslant i\leqslant K$, and let $$W_1=\prod_{i=1}^K
W_{1,i}.$$ Then
$$
W_1\leqslant 2^{\sum_{i=1}^K\sum_{j=1}^{h_i}p_{i,j}}\leqslant
2^{\frac{u^2}{\log u}},
$$
since (cf. \cite{SalatZnam68})
$$
\sum_{\substack{p\leqslant u\\ p\text{
prime}}}p=\bigg(\frac{1}{2}+o(1)\bigg)\frac{u^2}{\log u}.
$$
And noting that for sufficiently large $x$
$$
\frac{\log\log\log(2^{\frac{u^2}{\log u}})}{\log\log\log
(x^{\frac1K})}\leqslant\frac{2K(L+1)}{\log(\log\log x-\log
K)}\leqslant 1,
$$
we have $W_1\leqslant x^{\frac1K}$.

Let $m=\lfloor\log_2\log_2(x^{\frac2{K-1}})\rfloor$ and $
K'=1+\lfloor2^{-m}\log_2x\rfloor$, where $\log_2 x=\log x/\log 2$.
We have
$$
K'\leqslant1+\frac{\log_2x}{2^{m}}\leqslant1+\frac{2\log_2x}{2^{\log_2\log_2
(x^{\frac2{K-1}})}}=1+\frac{2\log_2x}{\frac2{K-1}\cdot\log_2 x}=K.
$$
For each $k\geqslant 0$, let $\gamma_k$ be the smallest prime factor
of $2^{2^k}+1$. Let
$$
W_2=\prod_{k=0}^{m-1}\gamma_k
$$
and $W=W_1W_2$. It is not difficult to see that $(W_1,W_2)=1$. And
$$
W\leqslant W_1\prod_{k=0}^{m-1}(1+2^{2^k})\leqslant x^{\frac1K}\cdot
x^{\frac2{K-1}}\leqslant x^{\frac3{K-1}}.
$$
Let $\beta$ be an odd integer such that
$$
\beta\equiv 2^{2^{m}(i-1)}+1\pmod{\prod_{j=1}^{h_i}q_{i,j}}
$$
and
$$
\beta\equiv 0\pmod{\gamma_k}
$$
for $1\leqslant i\leqslant K'$ and $0\leqslant k\leqslant m-1$.

Let
$$
\S=\{1\leqslant n\leqslant x:\, n\equiv\beta\pmod{2W}\}.
$$
Clearly, $$ \frac{x}{2W}-1\leqslant|\S|\leqslant \frac{x}{2W}. $$
Let
$$
\T_1=\{n\in\S:\, n\text{ is of the form } p+2^a+2^b\text{ with
}p\mid W\}
$$
and
$$
\T_2=\{n\in\S\setminus\T_1:\, n\text{ is of the form }
p+2^a+2^b\text{ with }p\nmid W\}.
$$
Clearly $|\T_1|=O(W(\log x)^2)$.

Suppose that $n\in\S$ and $n=p+2^a+2^b$ with $p$ is prime and
$0\leqslant a\leqslant b$. If $a\not\equiv b\pmod{2^m}$, then
$b=a+2^st$ where $0\leqslant s\leqslant m-1$ and $2\nmid t$. Thus
$$
p=n-2^a(2^{2^st}+1)\equiv\beta-2^a(2^{2^s}+1)\sum_{j=0}^{t-1}(-1)^j2^{2^sj}\equiv
0\pmod{\gamma_s}.
$$
Since $p$ is prime, we must have $p=\gamma_s$, i.e., $n\in\T_1$.

Below we assume that $a\equiv b\pmod{2^m}$. Write $b-a=2^m(t-1)$
where $1\leqslant t\leqslant K'$. If $a\equiv 0\pmod{p_{t,j}}$ for
some $1\leqslant j\leqslant h_t$, then recalling
$2^{p_{t,j}}\equiv1\pmod{q_{t,j}}$, we have
$$
p=n-2^a(2^{2^m(t-1)}+1)\equiv\beta-(2^{2^m(t-1)}+1)\equiv
0\pmod{q_{t,j}}.
$$
So $p=q_{t,j}$ and $n\in\T_1$. On the other hand, for any
$a\geqslant 0$ satisfying $a\not\equiv 0\pmod{p_{t,j}}$ for all
$1\leqslant j\leqslant h_t$, i.e., $(a,W_{1,t})=1$, by Lemma
\ref{l1}, we have
\begin{align*}
&|\{n\in\S:\,n-2^a(2^{2^m(t-1)}+1)\text{ is
prime}\}|\\\leqslant&\frac{2C_1|\S|}{\log
|\S|}\prod_{k=0}^{m-1}\bigg(1-\frac{1}{\gamma_k}\bigg)^{-1}\prod_{i=1}^K\prod_{j=1}^{h_i}\bigg(1-\frac{1}{q_{i,j}}\bigg)^{-1}\leqslant\frac{2^{5}C_1e^{2C_3}}{W}\cdot\frac{x}{\log
x}
\end{align*}
since $\gamma_k\equiv 1\pmod{2^{k+1}}$ and $\gamma_k>2^{k+1}$. And
noting that
$$
\frac{\log\log
u}{\log\log((\log_2x)^{\frac{1}{8}})}\leqslant\frac{K(L+1)}{\log(\log\log
x-\log\log2-\log 8)}<1,
$$
we have $u<(\log_2x)^{\frac{1}8}$. By Lemma \ref{l2},
\begin{align*}
&|\{0\leqslant a\leqslant\log_2x:\,a\not\equiv0\pmod{p_{t,j}}\text{
for all }1\leqslant j\leqslant h_t\}|\\\leqslant&C_2\frac{\log
x}{\log 2}\prod_{j=1}^{h_t}\bigg(1-\frac{1}{p_{t,j}}\bigg)\leqslant
2C_2e^{-L}\log x.
\end{align*}
Thus
\begin{align*}
|\T_2|\leqslant&\sum_{t=1}^{K'}\sum_{\substack{0\leqslant a\leqslant\log_2 x\\
(a,W_{1,t})=1}}|\{n\in\S:\,n-2^a(2^{2^m(t-1)}+1)\text{ is prime}\}|
\\\leqslant&K\cdot\frac{2^{5}C_1e^{2C_3}}{W}\cdot\frac{x}{\log x}\cdot2C_2e^{-L}\log x\leqslant\frac{x}{4W}.
\end{align*}
It follows that
\begin{align*}
&|\{n\in\S:\,n\text{ is not of the form
}p+2^a+2^b\}|\\=&|\S|-|\T_1|-|\T_2|\geqslant\frac{x}{2W}-1-O(W(\log
x)^2)-\frac{x}{4W}\gg x^{1-\frac{4}{K}}.
\end{align*}
The proof of Theorem \ref{t1} is complete. \qed

\begin{Rem}
Using a similar discussion, it is not difficult to prove that for
any given ${\mathcal {K}}\geqslant 1$,
\begin{align*}
&|\{1\leqslant n\leqslant x:\,n\text{ is odd and
}n\not=p+c(2^a+2^b)\text{ with }p\text{ prime},\ a,b\geqslant 0,\
1\leqslant c\leqslant
\mathcal {K}\}|\\
\gg&_{\mathcal {K}} x\cdot\exp\bigg(-C_{\mathcal {K}}\log
x\cdot\frac{\log\log\log\log x}{\log\log\log x}\bigg),
\end{align*}
where the constant $C_{\mathcal {K}}>0$ only depends on ${\mathcal
{K}}$.
\end{Rem}

\begin{Ack} The author is grateful to Professor Zhi-Wei Sun for his
helpful discussions on this paper.
\end{Ack}

\end{document}